# Asymptotic of summation arithmetic functions, the limit for which is the law of normal distribution

VICTOR VOLFSON

ABSTRACT

Summation arithmetic functions with asymptotically independent terms are studied in the paper, the limit of which is the law of normal distribution. Statements about the asymptotic behavior of the indicated functions are proved.

1. INTRODUCTION

In general, an arithmetic function is a function defined on the set of natural numbers and taking values on the set of complex numbers - $f: N \to C$. The name arithmetic function is connected with the fact that this function expresses a certain arithmetic property of the natural series.

Summation arithmetic function is a function:

$$S(x) = \sum_{n \leq x} f(n). \tag{1.1}$$

An example of summation arithmetic functions (1.1) is Mertens function:

$$M(x) = \sum_{k \leq x} \mu(k), \tag{1.2}$$

where $\mu(k)$ is Möbius function.

Möbius function $\mu(k) = 1$, if the positive integer k has an even number of prime divisors of the first degree, $\mu(k) = -1$, if the positive integer k has an odd number of prime divisors of the first degree and $\mu(k) = 0$, if the positive integer k has simple divisors not only of the first degree.

Keywords: summation arithmetic function, asymptotically independent terms, asymptotic behavior, Mertens function, Liouville function, Möbius function, limit distribution of summation arithmetic functions, law of normal distribution.

___

The study of summation arithmetic functions is usually carried out in two directions: the study of their distribution depending on the values of the natural argument and the study of their asymptotic behavior with the value of the natural argument $n \to \infty$.

$S(x) = \sum_{n \leq x} f(n)$ is a step function, which undergoes jumps at natural value of the argument, therefore with a certain accuracy it can be represented by a certain integral on the interval from 1 to $n$.

We can determine the asymptotic behavior $S(x) = \sum_{n \leq x} f(n)$ using the Euler-Maclaurin formula [1], [2] if $f$ is an elementary function.

If $f$ is not an elementary function, then the complex integration method is often used to determine the asymptotic behavior of $S(x) = \sum_{n \leq x} f(n)$ [3].

The study of the distribution of arithmetic functions is usually reduced to the determination of their average values [4].

The probabilistic approach to the study of the distribution of arithmetic functions [5] makes it possible to determine its variance, other moments of higher orders, the distribution function, and the characteristic function.

In this case, probabilistic spaces are determined $(\Omega_n, \mathcal{A}_n, \mathbb{P}_n)$ by taking $\Omega_n = \{1, 2, ..., n\}$, $\mathcal{A}_n$ - all the subsets $\Omega_n$, $P_n(A) = \frac{1}{n}\{N(m \in A\}$, where $N(m \in A)$ is the number of members of the natural series that satisfy the condition $m \in A$.

Then an arbitrary (real) arithmetic function $f$ (or rather, its restriction to $\Omega_n$) can be considered as a sequence of random variables $\xi_n$ defined on different probability spaces: $\xi_n(k) = f(k), 1 \leqslant k \leqslant n$.

When certain conditions are met, the specified sequence of random variables with the value $n \to \infty$ weakly converges in distribution to a certain distribution function. We say (in this case) that the limit for the arithmetic function is the given distribution function.

In [6], a theorem was proved that the limit distribution is an arithmetic function of the number of prime divisors of a natural number $n$ when value $n \to \infty$ is the function of the standard normal distribution. In [5], [7], [8], conditions were found under which the limit distribution of additive and multiplicative functions is the function of the standard normal distribution.

Definition. Asymptotic independence of arithmetic functions will be understood that the limit of the difference between the average value of the product of an arithmetic function $f(k)$ and the product of the average values of the same function for different values of the argument tends to zero when the value $n \to \infty$.

In [9], it was proved that the property of asymptotic independence of the terms holds for the summation functions of Mertens $M(x) = \sum_{n \leq x} \mu(n)$ and Liouville $L(x) = \sum_{n \leq x} \lambda(n)$. In [10], several assertions on classes of summation arithmetic functions whose terms have the property of asymptotic independence were proved. In addition (in [10]) conditions were found under which the normal distribution law for summation arithmetic functions is found to be the limit.

We also continue the study of summation arithmetic functions with asymptotically independent terms that have the limit normal distribution law, and define the asymptotic behavior of these functions.

The aim of the paper is to determine the existence of the limit normal distribution of the summation arithmetic function from its asymptotic.

2. ASYMPTOTIC BEHAVIOR OF SUMMATION ARITHMETIC FUNCTION, THE LIMIT FOR WHICH IS THE LAW OF NORMAL DISTRIBUTION

The assertion is proved in [1].

Assertion 1

Suppose there is a summation arithmetic function $S(n) = \sum_{k=1}^{n} f(k)$, where the term - arithmetic function $f : N \to R$ is bounded. Consider a sequence of random variables $f_n : f_n(k) = f(k), (n = 1, 2, ...)$ and a limiting random variable $f_0 : f_n \to f_0$ (by distribution when the value $n \to \infty$).

Suppose that for the mathematical expectation $f_n$ is performed:

$$M[f_n] = M[f_0] + o(1/n). \tag{2.1}$$

Then the limit for a sequence of random variables $S_n : S_n(k) = S(k), (1 \leq k \leq n)$ (or summation arithmetic function $S(n)$) is the law of normal distribution.

Assertion 2

A necessary and sufficient condition for the fulfillment of assertion 1 is that the asymptotic behavior of the sum of arithmetic functions $S(n) = \sum_{k=1}^{n} f(k)$ is:

$$S(n) = S_0(n) + o(1), \tag{2.2}$$

where $S_0(n) = nM[f_o]$, and $f(n)$ is limited arithmetic function.

Proof

Let assertion 1 be fulfilled, then for a sequence of random variables $f_n$ and a random variable $f_o : f_n \to f_o$ (by distribution when the value $n \to \infty$) (2.1) is fulfilled.

Considering what $M[f_n] = \sum_{k=1}^{n} f_n(k)/n$ and $f_n(k) = f(k), (1 \leq k \leq n)$ we get:

$$M[f_n] = \sum_{k=1}^{n} f(k) = M[f, n], \tag{2.3}$$

where $M[f, n]$ is the average value of the arithmetic function $f$ on the interval from 1 to $n$.

Since $M[f_n] = \sum_{k=1}^{n} f_n(k)/n = S(n)/n$ and having in mind (2.1):

$$S(n)/n = M[f_0] + o(1/n),$$

therefore:

$$S(n) = nM[f_o] + o(1) = S_o(n) + o(1),$$

which corresponds to (2.2).

The condition of boundedness of the arithmetic function $f$ is also fulfilled based on statement 1.

We now prove that if (2.2) holds and $f(n)$ is a bounded arithmetic function, then the conditions of assertion 1 are satisfied.

The condition of boundedness of the arithmetic function $f(n)$ in assertion 1 is satisfied.

Let's divide in the formula (2.2) the right and left parts on $n$ and we get:

$$S(n)/n = \sum_{k=1}^{n} f(k)/n = \sum_{k=1}^{n} f_n(k) = M[f_n].$$

On the other hand $S_0(n) = nM[f_o]$, therefore:

$$S(n)/n = S_0(n)/n + o(1/n) \text{ or } M[f_n] = M[f_o] + o(1/n).$$

Therefore, all the conditions of assertion 1 are satisfied.

Assertion 3

If the summand function $f$ in the summation arithmetic function $S(n) = \sum_{k=1}^{n} f(k)$ is bounded and elementary, then assertion 1 is not fulfilled.

Proof

If the function $f$ is bounded and elementary, then 2 cases are possible: $f$ is a decreasing function, $f$ is not a decreasing function.

If $f$ is a decreasing function, then based on [2]:

$$\sum_{k=1}^{n} f(k) = \int_{1}^{n} f(t)dt + O(1),$$

therefore formula (2.2) of assertion 2 and, accordingly, assertion 1 are not satisfied.

If $f$ is not a decreasing function, then based on [2]:

$$\sum_{k=1}^{n} f(k) = \int_{1}^{n} f(t)dt + O(f(n)). \qquad (2.4)$$

Since the function $f$ is bounded and based on $f(n) = O(1)$ then (2.4) can also be written as:

$$\sum_{k=1}^{n} f(k) = \int_{1}^{n} f(t)dt + O(1),$$

therefore, in this case, formula (2.2) of assertion 2 and, accordingly, assertion 1 are not satisfied.

Having in mind assertion 3, assertion 1 can be performed only when the function $f$ in the summation arithmetic function $S(n) = \sum_{k=1}^{n} f(k)$ is bounded, but it is not an elementary function. Let us consider this case.

Assertion 4

Let the summation arithmetic function $S(n) = \sum_{k=1}^{n} f(k)$ has $\lim_{n \to \infty} S(n) = 0$, and the arithmetic function $f : N \to R$ is bounded.

Then the limit for the summation arithmetic function $S(n)$ is the normal distribution law and $M[f_0] = 0$.

Proof

Based on the condition - $\lim_{n \to \infty} S(n) = 0$, therefore:

$$S(n) = o(1). \qquad (2.5)$$

Divide (2.5) by $n$ and get:

$$\frac{S(n)}{n} = M[f,n] = M[f_n] = o(1/n). \tag{2.6}$$

Thus, all the conditions of assertion 2 with $M[f_0]=0$ are satisfied. Therefore, assertion 1 is fulfilled for summation arithmetic function $S(n) = \sum_{k=1}^{n} f(k)$, therefore, the law of normal distribution is the limiting one for this summation function.

Two cases are possible, when performing assertion 4:

1. The arithmetic function $f: N \to R$ is decreasing and $\lim_{n \to \infty} f(n) = 0$.

2. The arithmetic function $f: N \to R$ is a bounded and non-decreasing function.

Let us consider the example of the first case. Based on the work [11] $\lim_{n \to \infty} \sum_{k=1}^{n} \frac{\mu(k)}{k} = 0$. On the other hand $|\frac{\mu(n)}{n}| \leq 1$, i.e. $f: N \to R$ is limited. Therefore, all conditions of assertion 2 are fulfilled.

Let us consider the example of the second case. Suppose there is a summation arithmetic function - $S(n) = \sum_{k=1}^{n} f(k)$, where the arithmetic function $f: N \to R$ corresponds to a sequence of random variables $f_n (n = 1, 2, ...)$ (in the probability spaces considered in Section 1).

Let a random variable $f_n = 1$ with probability - $p_{1n} = 1/2 + 1/(n+1)\log(n+1)$ and $f_n = -1$ with probability - $p_{2n} = 1/2 - 1/(n+1)\log(n+1)$ ($p_{1n} + p_{2n} = 1$).

The mathematical expectation $f_n$ is:

$$M[f_n] = 1/2 + 1/(n+1)\log(n+1) - 1/2 + 1/(n+1)\log(n+1) = 2/(n+1)\log(n+1). \tag{2.7}$$

Having in mind that $M[f_n] = S(n)/n$, based on (2.7) we get:

$$S(n) = 2/\log(n+1) \text{ and } \lim_{n \to \infty} S(n) = \lim_{n \to \infty} 2/\log(n+1) = 0.$$

Since $f: N \to R$ is a bounded arithmetic function, all the conditions of Statement 4 are satisfied.

Let us generalize the second case, when the arithmetic function $f: N \to R$ is bounded non-decreasing, and the limit for the summation arithmetic function $S(n) = \sum_{k=1}^{n} f(k)$ is the law of normal distribution.

Assertion 5

Suppose there is a summation arithmetic function - $S(n) = \sum_{k=1}^{n} f(k)$, where the arithmetic function $f: N \to R$ corresponds to a sequence of random variables $f_n (n=1,2,...)$ (in the probability spaces considered in Section 1).

Let us a random variable $f_n$ take values $a_1,...,a_l$ with probabilities, respectively: $p_{1n} = p_{10} + o(1/n),..., p_{ln} = p_{l0} + o(1/n)$. Then the limit for summation arithmetic functions $S(n) = \sum_{k=1}^{n} f(k)$ is the law of normal distribution, and the asymptotic of summation arithmetic functions will be equal to:

$$S(n) = n(\sum_{k=1}^{n} a_i p_{io}) + o(1). \qquad (2.8)$$

Proof

Based on the condition, the function $f: N \to R$ is limited. Thus, the first condition of assertion 1 is satisfied. Let us prove the fulfillment of the second condition of assertion 1:

$$M[f_n] = M[f,n] = \sum_{i=1}^{l} a_i [p_{io} + o(1/n)] = \sum_{i=1}^{l} a_i p_{i0} + o(1/n) = M[f_o] + o(1/n), \qquad (2.9)$$

where $M[f_o] = \sum_{i=1}^{l} a_i p_{io}$.

Thus, all the conditions of assertion1 are satisfied; therefore, the summation arithmetic function $S(n) = \sum_{k=1}^{n} f(k)$ has a limit normal distribution law.

When assertion 1 is fulfilled, the summation arithmetic function has an asymptotic behavior in accordance with assertion 2 and having in mind (2.9):

$$S(n) = nM[f_o] + o(1) = n(\sum_{i=1}^{l} a_i p_{io}) + o(1),$$

which corresponds to (2.8).

Let us consider an example of the fulfillment of assertion 5, but in contrast to the second example $M[f_o] \neq 0$.

Suppose there is summation arithmetic function - $S(n) = \sum_{k=1}^{n} f(k)$, where the arithmetic function $f : N \to R$ corresponds to a sequence of random variables $f_n (n = 1, 2, ...)$ (in the probability spaces considered in Section 1).

Пусть случайная величина $f_n = 1$ с вероятностью $p_{1n} = 1/2 + 1/(n+1)\log^2(n+1)$ и $f_n = 0$ с вероятностью $p_{2n} = 1/2 - 1/(n+1)\log^2(n+1)$.

Let a random variable $f_n = 1$ with probability $p_{1n} = 1/2 + 1/(n+1)\log^2(n+1)$ and $f_n = 0$ with probability $p_{2n} = 1/2 - 1/(n+1)\log^2(n+1)$.

Следовательно, $f : N \to R$ удовлетворяет всем условия утверждения 5 и предельным для сумматорной арифметической $S(n) = \sum_{k=1}^{n} f(k)$ является закон нормального распределения, и асимптотика сумматорной арифметической функции в этом случае будет: $S(n) = n/2 + O(1/\log^2(n))$.

Therefore, $f : N \to R$ satisfies all the conditions of assertion 5 and the limit for the summation arithmetic function $S(n) = \sum_{k=1}^{n} f(k)$ is the law of normal distribution, and the asymptotic of the summation arithmetic function in this case will be: $S(n) = n/2 + O(1/\log^2(n))$.

3. CONCLUSION AND SUGGESTIONS FOR FURTHER WORK

The next article will continue to study the behavior of some sums.

## 4. ACKNOWLEDGEMENTS

Thanks to everyone who has contributed to the discussion of this paper. I am grateful to everyone who expressed their suggestions and comments in the course of this work.